\documentclass{amsart}
\usepackage{amssymb,latexsym}

\title{Stokes Theorem for Lipschitz forms on a smooth manifold}

\author{Stanislav Dubrovskiy}

\address{Department of Mathematics, Ben Gurion
University of the Negev, P.O.B. 653\\Beer'Sheva 84105, Israel}

\email{dubr@math.bgu.ac.il}

\keywords{Lipschitz forms, Stokes theorem}

\date{May 27, 2008}

\begin{document}

\maketitle

\begin{abstract}
Stokes theorem holds for Lipschitz forms on a smooth manifold.
\end{abstract}

Let $\omega$ be a Lipschitz $(n-1)$-form with compact support $K$
on a smooth\\ $n$-manifold $M$, in the notation of \cite{GKS}:
$$\omega \in \mathcal{L}^{n-1}(M)\ .$$

Suppose first $M=\mathbb{R}^{n}$ or $\mathbb{H}^{n}=\{x_n\geq
0\}$, the upper half space in $\mathbb{R}^{n}$.\\It suffices to
assume $$\omega=fdg_1\wedge dg_2\ldots\wedge dg_{n-1}$$ (rather
than a finite sum of such terms), where $f$, $g_1$, $g_2,\ldots$,
$g_{n-1}$ are\\ Lipschitz functions on $M$. According to
\cite[Lemma 1.1]{GKS}:
$$d\omega=df\wedge dg_1\wedge dg_2\ldots\wedge dg_{n-1}\ .$$
We can write $$\omega=a\,dx_1\wedge dx_2\ldots\wedge
dx_{n-1}+\omega_n\ ,$$ where $\omega_n$ contains $dx_n$, and $a$
is a polynomial in $f$, $\frac{\partial g_i}{\partial x_j}$.\\Then
$$d\omega=b\,dx_1\wedge dx_2\ldots\wedge dx_n\ ,$$ $b$ is a polynomial
in $\frac{\partial f_i}{\partial x_j}$, $\frac{\partial
g_i}{\partial x_j}$.

Consider $a|_{x_n=0}$, the restriction of $a$ to $K'=K\cap
\{x_n=0\}$, a compact subset in $K$. It is again a polynomial in
Lipschitz functions, since restrictions of Lipschitz functions to
a hyperplane are Lipschitz.

Just as in the proof of \cite[Lemma 1.1]{GKS} let us consider
smooth forms
$$\omega_s=f_{\varepsilon_s}d(g_1)_{\varepsilon_s}\wedge
d(g_2)_{\varepsilon_s}\ldots\wedge d(g_{n-1})_{\varepsilon_s}\
$$and$$\ d\omega_s=df_{\varepsilon_s}d(g_1)_{\varepsilon_s}\wedge
d(g_2)_{\varepsilon_s}\ldots\wedge d(g_{n-1})_{\varepsilon_s}$$ -
sequence $\varepsilon_s$ chosen so that the Steklov averages
$f_{\varepsilon_s}$, $(g_1)_{\varepsilon_s}$,
$(g_2)_{\varepsilon_s},\ldots$, $(g_{n-1})_{\varepsilon_s}$ have
their derivatives converge:
$$\left(\frac{\partial f}{\partial x_j}\right)_{\varepsilon_s}\rightarrow\frac{\partial f}{\partial x_j}\
,\quad \left(\frac{\partial g_i}{\partial
x_j}\right)_{\varepsilon_s}\rightarrow\frac{\partial g_i}{\partial
x_j}$$ almost everywhere on $K$ and $K'$. \\(This can be done
since the restrictions to $K'$ are still Lipschitz functions.)\\

Functions $$a(f_{\varepsilon_s},(\,\frac{\partial g_i}{\partial
x_j}\,)_{\varepsilon_s})\quad ,\quad b\,(\,(\,\frac{\partial
f_i}{\partial x_j}\,)_{\varepsilon_s},(\,\frac{\partial
g_i}{\partial x_j}\,)_{\varepsilon_s})$$ are uniformly bounded in
$s$ and converge almost everywhere to $$a(f,\frac{\partial
g_i}{\partial x_j})\ \quad ,\quad b\,(\,\frac{\partial
f_i}{\partial x_j}\,,\frac{\partial g_i}{\partial x_j}\,)$$
respectively.

Then by the Lebesgue limit theorem:
$$\int_{\partial M} \omega_s \rightarrow \int_{\partial M} \omega \quad, \quad \int_{M} d(\omega_s) \rightarrow \int_{M} d\omega
\ .$$ Since Stokes theorem is valid for smooth $\omega_s$, we must
conclude it holds for the Lipschitz forms in the limit.

Suppose now $M$ consists of a single chart $\phi: M\rightarrow
\mathbb{R}^{n} (\textrm{or}\ \mathbb{H}^{n})$. \\We have:
$$\int_{\partial M} \omega = \int_{\partial (\phi(M))}
(\phi^{-1})^{*}\omega = \int_{\phi(M)} (\phi^{-1})^{*}d\omega =
\int_{M} d\omega \ .
$$
For an atlas with multiple charts, we use partition of unity
$\{\rho_\alpha\}$ and the fact that $\rho_\alpha \omega$ are again
Lipschitz, together with the linearity of integral, cf. \cite{BT}.

Integration on Lipschitz manifolds is to be discussed in a
forthcoming paper.

\end{document}